\documentclass[12pt]{article}
\usepackage{amsthm,amsfonts,amssymb,amsmath}
\usepackage{stmaryrd}
\usepackage{cite,hyperref}
\usepackage{epsfig}
\usepackage{url}
\usepackage{xcolor,tikz}
\usetikzlibrary{positioning, arrows.meta,calc}
\usetikzlibrary{matrix}
\usepackage{multicol,graphicx}
\usepackage{fullpage}
\usepackage{bbm}
\usetikzlibrary{decorations}
\usepackage{svg}

\usepackage[shortlabels]{enumitem}

\numberwithin{equation}{section}

\newtheorem{thm}[equation]{Theorem}

\theoremstyle{definition}

\newtheorem*{ai}{Acknowledgements and AI Declaration}

\theoremstyle{remark}

\title{The Partial List Colouring Conjecture is False}
\author{Jonathan A. Noel\thanks{Department of Mathematics and Statistics, University of Victoria, Victoria, B.C., Canada. E-mail: {\tt noelj@uvic.ca}. Research supported by NSERC Discovery Grant RGPIN-2021-02460.}}

\DeclareTextCompositeCommand{\v}{OT1}{l}{l\nobreak\hspace{-.1em}'}
\DeclareTextCompositeCommand{\v}{OT1}{t}{t\nobreak\hspace{-.1em}'\nobreak\hspace{-.15em}}

\begin{document}

\maketitle

\begin{abstract}
We exhibit a graph $G$ with $14$ vertices and list chromatic number equal to $3$ such that there is a $2$-list assignment $L$ of $G$ such that at most $9$ vertices of $G$ can be properly coloured from $L$. This disproves the Partial List Colouring Conjecture of Albertson, Grossman and Haas. This counterexample was discovered and fully verified by ChatGPT 6 Astra Ultra after some persistent prompting, but almost no mathematical input, from the author.
\end{abstract}

\section{Introduction}

Let $G$ be a graph with $n$ vertices and chromatic number $k\geq 2$ and let $1\leq \ell< k$. Suppose that your goal is to properly colour as many vertices of $G$ as possible using only $\ell$ colours. Trivially, you can colour at least $n(\ell/k)$ vertices by simply taking the largest $\ell$ colour classes of a proper $k$-colouring. Albertson, Grossman and Haas~\cite{AlbertsonGrossmanHaas00} conjectured that a similar phenomenon is true for list colouring. That is, if $G$ has $n$ vertices and list chromatic number equal to $k\geq2$ and $1\leq \ell<k$, then, for any $\ell$-list assignment $L$ of $G$, it is possible to properly $L$-colour at least $n(\ell/k)$ vertices. This has come to be known as the \emph{Partial List Colouring Conjecture}.

There are several partial results on this conjecture. Chappell~\cite{Chappell99} and Haas, Hanson and MacGillivray~\cite{HaasHansonMacGillivray03} proved that one can always properly $L$-colour at least $\frac{6}{7}n(\ell/k)$ vertices and $n/\lceil k/\ell\rceil$ vertices, respectively. In particular, the result of~\cite{HaasHansonMacGillivray03} implies the conjecture in the case that $\ell\mid k$. Janssen~\cite{Janssen01} proved that the conjecture holds if $\Delta(G)\leq k$ or if $\left|\bigcup_{v\in V(G)}L(v)\right|\leq k$ and Janssen, Mathew and Rajendraprasad~\cite{JanssenMathewRajendraprasad15} proved it for various other classes of graphs. Iradmusa~\cite{Iradmusa10} proved that, for any $G$, the conjecture holds for at least half of the values $\ell\in\{1,\dots,k-1\}$. In spite of the above evidence, the conjecture remained open even in the first non-trivial case $k=3$ and $\ell=2$. Here, we prove that the conjecture is false in this case. 

\begin{thm}
\label{th:main}
There exists a graph $G$ with $14$ vertices and list chromatic number equal to $3$ and a $2$-list assignment $L$ of $G$ such that at most $9$ vertices of $G$ can be properly coloured from $L$.
\end{thm}

This disproves the conjecture of~\cite{AlbertsonGrossmanHaas00} since $9<\frac{28}{3}=14\left(\frac{2}{3}\right)$. In the next section, we describe the counterexample, verify its properties and make some concluding remarks.

\section{The Counterexample}

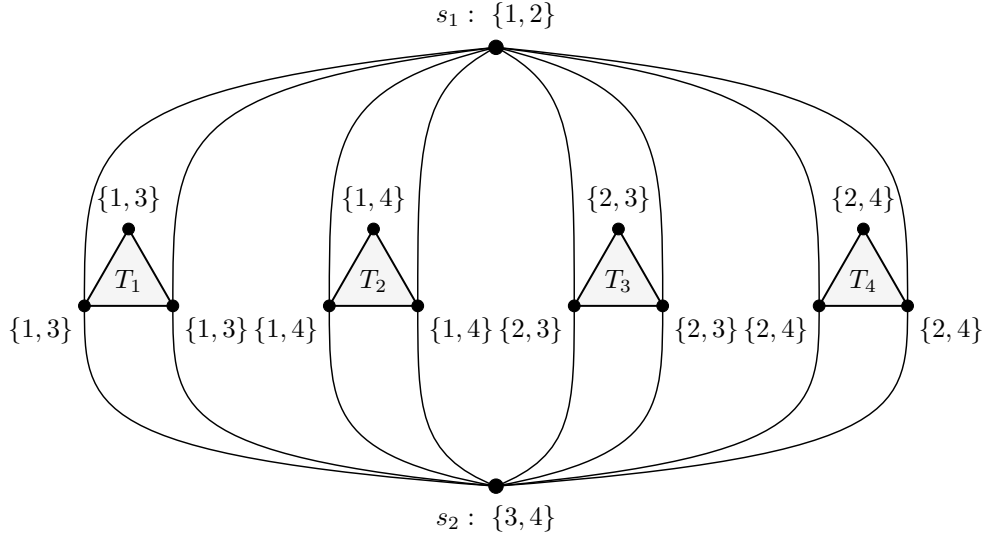
\begin{figure}[htbp]
    \centering
    \begin{tikzpicture}[scale=0.9,
  x=1cm,y=1cm,
  edge/.style={draw=black,line width=.55pt},
  clique edge/.style={draw=black,line width=.75pt},
  vertex/.style={circle,fill=black,draw=black,inner sep=0pt,
                 minimum size=4.4pt},
  terminal/.style={vertex,minimum size=5.4pt},
  list/.style={font=\footnotesize,inner sep=1pt},
  clique label/.style={font=\footnotesize,inner sep=0pt}
]
  \coordinate (t1) at (0,3.8);
  \coordinate (t2) at (0,-2.65);

  \foreach \i/\x/\a/\b in {
      1/-5.4/1/3,
      2/-1.8/1/4,
      3/1.8/2/3,
      4/5.4/2/4}{
    \coordinate (u\i) at ({\x-.65},0);
    \coordinate (v\i) at ({\x+.65},0);
    \coordinate (p\i) at (\x,1.13);

    \foreach \dx/\name in {-.65/u,.65/v}{
      \draw[edge] (t1)
        .. controls ({\x+\dx},3.2) and ({\x+\dx},2.4)
        .. (\name\i);
      \draw[edge] (t2)
        .. controls ({\x+\dx},-2.2) and ({\x+\dx},-1.5)
        .. (\name\i);
    }

    \path[fill=black!4] (u\i)--(v\i)--(p\i)--cycle;
    \draw[clique edge] (u\i)--(v\i)--(p\i)--cycle;
    \node[clique label] at (\x,.36) {$T_{\i}$};

    \node[vertex] at (u\i) {};
    \node[vertex] at (v\i) {};
    \node[vertex] at (p\i) {};
    \node[list,anchor=north east,xshift=-3pt,yshift=-3pt]
      at (u\i) {$\{\a,\b\}$};
    \node[list,anchor=north west,xshift=3pt,yshift=-3pt]
      at (v\i) {$\{\a,\b\}$};
    \node[list,above=5pt] at (p\i) {$\{\a,\b\}$};
  }

  \node[terminal] at (t1) {};
  \node[terminal] at (t2) {};
  \node[list,above=6pt] at (t1) {$s_1:\ \{1,2\}$};
  \node[list,below=6pt] at (t2) {$s_2:\ \{3,4\}$};
\end{tikzpicture}
    \caption{The counterexample and a $2$-list assignment from which at most $9$ vertices can be properly coloured.}
    \label{fig:graph}
\end{figure}

The graph $G$ is constructed as follows. Let $T_1,T_2,T_3,T_4$ be pairwise disjoint triangles and let $p_i$ be a vertex of $T_i$, which we call the \emph{private} vertex of $T_i$, for each $1\leq i\leq 4$. Let $s_1$ and $s_2$ be two new non-adjacent vertices that we call \emph{terminals} and add all edges from the terminals to non-private vertices of $T_i$ for $1\leq i\leq 4$. Let $L$ be the $2$-list assignment in which 
\begin{itemize}
\item $L(s_1)=\{1,2\}$, 
\item $L(s_2)=\{3,4\}$, 
\item every vertex of $T_1$ has list $\{1,3\}$, 
\item every vertex of $T_2$ has list $\{1,4\}$, 
\item every vertex of $T_3$ has list $\{2,3\}$, and
\item every vertex of $T_4$ has list $\{2,4\}$.
\end{itemize}
In other words, each of the possible lists of size two that includes one colour from $L(s_1)$ and one from $L(s_2)$ is assigned to one of the triangles $T_1,T_2,T_3,T_4$. See Figure~\ref{fig:graph} for a depiction of the graph $G$ and this list assignment. 

First, let us argue that it is not possible to properly colour more than $9$ vertices of $G$ from the assignment $L$. For $1\leq i\leq 4$, all vertices of $T_i$ have the same list of size two. So, any partial proper colouring from the lists $L$ cannot colour more than two vertices from each triangle. Thus, in order to colour $10$ vertices, both terminals need to be coloured. However, any combination of colours of the terminals appears as the list assigned to one of the triangles. Thus, if both terminals are coloured, then there exists one of the triangles, say $T_i$, such that at most one vertex of $T_i$ is coloured. Since at most two vertices of the other three triangles are coloured, the total number of vertices coloured is at most $9$. 

Finally, we argue that $G$ has list chromatic number equal to three. Since it contains a triangle, its list chromatic number is at least three. Now, let $L$ be any $3$-list assignment for $G$. Suppose first that $L(s_1)\cap L(s_2)\neq\emptyset$. In this case, we simply colour $s_1$ and $s_2$ with a colour that is common to their lists. After doing so and deleting that colour from the lists of the neighbours of $s_1$ and $s_2$, each triangle still has one vertex with a list of size three and two others with lists of size at least two. So, the triangles can be coloured greedily. 

So, $L(s_1)\cap L(s_2)=\emptyset$. For $1\leq i\leq 4$, let $v_i$ be an arbitrary non-private vertex of $T_i$. Then, since $L(s_1)\cap L(s_2)=\emptyset$, we have $|L(v_i)\cap L(s_1)|\cdot |L(v_i)\cap L(s_2)|\leq 2$. Thus, the sum of this quantity over all $i\in \{1,2,3,4\}$ is at most $8$. Since $|L(s_1)|\cdot |L(s_2)|=9$, there must exist a pair $(c,d)\in L(s_1)\times L(s_2)$ such that there does not exist $1\leq i\leq 4$ with $c,d\in L(v_i)$. Colour $s_1$ with $c$ and $s_2$ with $d$ and remove these colours from the lists of all of their neighbours. After doing so, for each $1\leq i\leq 4$, the list of $p_i$ still has three elements, the list of $v_i$ still has at least two elements, and the list of the other vertex of $T_i$ still has at least one element. So, the triangles can be greedily coloured from their remaining lists. This completes the proof of Theorem~\ref{th:main}.

We remark that the additive gap between the number of properly $L$-colourable vertices in this counterexample and the number predicted by the Partial List Colouring Conjecture is very small. To inflate it, one can simply take a disjoint union of arbitrarily many copies of this construction. It may be interesting to try to cook up connected counterexamples which achieve a large gap from the conjecture.

\begin{ai}
The author has believed that the Partial List Colouring Conjecture is likely to be false since the late 2010s and has tried various human and computer searches to disprove it including, more recently, various AI-assisted tools. On September 18, 2026, he asked ChatGPT 6 Astra Ultra to disprove it. It was unable to do so initially but, after a few additional prompts and some persistence, it eventually found a counterexample which required extensive computer assistance to verify. The author then pushed it to find a human-checkable proof, and it found one that was moderately complicated. After some more pushing, it eventually found the intuitive construction with the simple proof that we have included here. This paper has been written from scratch based on a draft that was generated by AI, with no further AI assistance other than generating the picture in Figure~\ref{fig:graph} and a bit of proofreading. All arguments were fully checked by the author, who takes full responsibility for correctness. 
\end{ai}

\end{document}